\documentclass[12pt]{amsart}
\usepackage{amssymb,amsmath}
\numberwithin{equation}{section}
\def\sumj{\underset{j=1}{\overset n\sum}}

\def\di{\partial}
\def\simleq{\underset\sim<}
\def\simgeq{\underset\sim>}

\def\T{\text}
\def\1#1{\overline{#1}}
\def\2#1{\widetilde{#1}}
\def\3#1{\widehat{#1}}
\def\4#1{\mathbb{#1}}
\def\5#1{\frak{#1}}
\def\6#1{{\mathcal{#1}}}
\def\C{{\4C}}
\def\R{{\4R}}

\def\Re{{\sf Re}\,}
\def\Im{{\sf Im}\,}

\def\phi{\varphi}

\newtheorem{Thm}{Theorem}[section]
\newtheorem{Cor}[Thm]{Corollary}
\newtheorem{Pro}[Thm]{Proposition}
\newtheorem{Lem}[Thm]{Lemma}
\theoremstyle{definition}\newtheorem{Def}[Thm]{Definition}
\theoremstyle{remark}
\newtheorem{Rem}[Thm]{Remark}
\newtheorem{Exa}[Thm]{Example}

\def\Label#1{\label{#1}}
\def\bl{\begin{Lem}}
\def\el{\end{Lem}}
\def\bp{\begin{Pro}}
\def\ep{\end{Pro}}
\def\bt{\begin{Thm}}
\def\et{\end{Thm}}
\def\bc{\begin{Cor}}
\def\ec{\end{Cor}}
\def\bd{\begin{Def}}
\def\ed{\end{Def}}
\def\br{\begin{Rem}}
\def\er{\end{Rem}}
\def\be{\begin{Exa}}
\def\ee{\end{Exa}}
\def\bpf{\begin{proof}}
\def\epf{\end{proof}}
\def\ben{\begin{enumerate}}
\def\een{\end{enumerate}}

\def\1alpha{[\frac1\alpha]}
\def\T{\text}
\def\R{{\Bbb R}}

\def\C{{\Bbb C}}

\numberwithin{equation}{section}
\def\T{\text}

\newcommand{\no}[1]{\|{#1}\|}
\newcommand{\NO}[1]{{\|{#1}\|}^2}
\newtheorem{theorem}{Theorem  }[section]
\newtheorem{definition}[theorem]{Definition }
\newtheorem{lemma}[theorem]{Lemma  }
\newtheorem{proposition}[theorem]{Proposition  }
\newtheorem{corollary}[theorem]{Corollary }
\newtheorem{example}[theorem]{\it Example }

\begin{document}
\title[Loss of derivatives...]{Loss of derivatives for systems of complex vector fields and sums of squares}        
\author[ T.V.~Khanh, S.~Pinton and G.~Zampieri ]
{Tran Vu Khanh, Stefano Pinton and Giuseppe Zampieri}
\address{Dipartimento di Matematica, Universit\`a di Padova, via 
Trieste 63, 35121 Padova, Italy}
\email{khanh@math.unipd.it, pinton@math.unipd.it, 
zampieri@math.unipd.it}
\maketitle

\begin{abstract}
We discuss, both for systems of complex vector fields and for sums of squares, the phenomenon discovered by Kohn of hypoellipticity with loss of derivatives.
\vskip0.2cm
\noindent
MSC: 32W05, 32W25, 32T25

\end{abstract}
\section{Estimates for vector fields and sums of squares in $\R^3$}
\Label{s1}
A system of real vector fields $\{X_j\}$ in $T\R^n$ is said to satisfy the bracket finite type condition if
\begin{equation}
\Label{bracket}
\T{commutators of order $\leq h-1$ of the $X_j$'s span the whole $T\R^n$.}
\end{equation}
Explicitly: $\T{Span}\{X_j,\,[X_{j_1},X_{j_2}],...,[X_{j_1},[X_{j_2},...,[X_{j_{h-1}},X_{j_h}]]...]\}=T\R^n$.
This system enjoys $\delta$-subelliptic estimates for $\delta=\frac1h$ and therefore it is hypoelliptic  according to H\"ormander \cite{H67}. 
(See also \cite{FK72} and \cite{KN65} for elliptic regularization which yields regularity from estimates.)
This remains true for systems of complex vector fields $\{L_j\}$ stable under conjugation (both in $\C\otimes T\R^n$ or $\C\otimes T\C^n$) once one applies H\"ormander's result to $\{\Re L_j,\,\Im L_j\}$. Stability under conjugation can be artificially achieved by adding $\{\epsilon\bar L_j\}$ in order to apply H\"ormander's theorem $\NO{u}_{\delta}\leq \sum_j(c_\epsilon\NO{L_ju}+\epsilon\NO{\bar L_ju})+c_\epsilon\NO{u},\,\,u\in C^\infty_c$. 
(Precision about $\epsilon$ and $c_\epsilon$ is not in the statement but transparent from the proof.)
On the other hand, by integration by parts $\NO{\bar L_ju}\simleq \NO{L_ju}+|([L_j,\bar L_j]u,u)|+\NO{u}\simleq \NO{L_ju}+\NO{u}_{\frac12}+ \NO{u}$. Thus if the type is $h=2$, and hence $\delta=\frac12$, the $\frac12$-norm is abbsorbed in the left:  $\{\epsilon\bar L_j\}$ can be taken back and one has $\frac12$-subelliptic estimates for $\{L_j\}$.
 The restraint $h=2$ is substantial and in 
fact Kohn discovered in \cite{K05} a pair of vector fields $\{L_1,\,L_2\}$ in $\R^3$ of finite type $k+1$ (any fixed $k$) which are not subelliptic but, 
nonetheless, are hypoelliptic. Precisely, in the terminology of \cite{K05}, they loose $\frac{k-1}2$ derivatives and the related sum of squares $\bar L_1L_
1+\bar L_2L_2$ looses $k-1$ derivatives. The vector fields in question are $L_1=\di_{\bar z}+iz\di_t$ and $L_2= \bar z^k(\di_z-i\bar z\di_t)$ in $\C\times \R$. 
Writing $t=\Im w$, they are identified to $\bar L$ and $\bar z^kL$ for the CR vector field $\bar L$ tangential to the strictly pseudoconvex hypersurface $
\Re w=|z|^2$ of $\C^2$. Consider a more general hypersurface $M\subset\C^2$ defined by $\Re w=g(z)$ for $g$ real, and use the notations $g_1=\di_zg$, $ g_{1
\bar1}=\di_{z}\di_{\bar z}g$ and $g_{1\bar1\bar1}=\di_z\di_{\bar z}\di_{\bar z}g$. Suppose  that $M$ is pseudoconvex, that is,  $g_{1\bar1}\geq0$ and 
denote by $2m$ the vanishing order of $g$ at $0$, that is, $g=0^{2m}$. Going further in the analysis of loss of derivatives, Bove, Derridj, Kohn and Tartakoff have 
considered the case where
\begin{equation}
\Label{b}
g_1=\bar z|z|^{2(m-1)}h(z) \T{ and } g_{1\bar1}=|z|^{2(m-1)}f(z)\T{ for  $f>0$.}
\end{equation}
If $L=\di_z-ig_1\di_t$ is the $(1,0)$ vector field tangential to $Re w=g$ for $g$ satisfying \eqref{b}, they have proved loss of $\frac{k-1}m$ derivatives for the operator $L\bar L+\bar L|z|^{2k}L$. 

We consider here a general pseudoconvex hypersurface  $M \subset\C^2$;   $\zeta$ and $\zeta'$ will denote cut-off functions in a neigborhood of $0$ such that $\zeta'|_{\T{supp}\,\zeta}\equiv1$.
\bt
\Label{p1.1}
Let $\{L,\,\bar L\}$ (or better $\{\Re L,\,\Im L\}$) have type $2m$; then the system $\{\bar L,\,\bar z^kL\}$ looses $l:=\frac{k-1}{2m}$ derivatives. More precisely
\begin{equation}
\Label{1.2}
\begin{split}
\no{\zeta u}^2_{s}&\simleq \no{\zeta'\bar L  u}^2_{s-\frac1{2m}}+\no{\zeta' \bar z^{k}\bar L  u}^2_{s+l  }
\\
&+\no{\zeta' \bar z^{k}L  u}^2_{s+l }+\no{u}^2_0.
\end{split}
\end{equation}
\et
The estimate \eqref{1.2} says that the responsible of the loss $l$ is $\bar z^kL$ (plus the extra vector field $\bar z^k\bar L$) and not $\bar L$. The proof of this here, as well as the two theorems below, follows in Section~\ref{s4}. What underlies the whole technicality is the basic notion of subelliptic multiplier; also the stability of multipliers under radicals is crucial (hidden in the interpolation Lemma~\ref{l2.1} below). We point out that though the coefficient of the vector field $\bar L$ gains much in generality ($+ig_{\bar1}$ instead of $+i z$ or $+i z|z|^{2(m-1)}$ as in \cite{K05} and \cite{BDKT06} respectively), instead, the perturbation $\bar z^k$ of $L$ remains the same. This is substantial; only an antiholomorphic  perturbation is allowed.
We introduce a new notation for the perturbed Kohn-Laplacian
\begin{equation}
\Label{1.3}
\Box^k=L\bar L+\bar L|z|^{2k}L\qquad\T{for $L=\di_z-ig_1\di_t$}.
\end{equation}
\bt
\Label{t1.1}
Let $\{L,\,\bar L\}$ have type $2m$ and assume moreover, that
\begin{equation}
\Label{1.4}
|g_1|\simleq |z|g_{1\bar1}\,\,\,\T{ and }\,\,\,|g_{1\bar1\bar1}|\simleq |z|^{-1}g_{1\bar1}.
\end{equation}
Then $\Box^k$ looses $l=\frac{k-1}m$ derivatives, that is
\begin{equation}
\Label{1.5}
\no{\zeta u}^2_s\simleq \no{\zeta'\Box^ku}^2_{s+2l}+\no{u}^2_0.
\end{equation}
\et
Differently from vector fields, loss for sums of squares requires the additional assumption \eqref{1.4}; whether finite type suffices is an open question.
\be
Consider the boundary defined by $\Re w=g$ with $g(z)=0^{2m}$ and assume
\begin{equation}
\Label{a}
g_{1\bar 1}\simgeq|z|^{2(m-1)}.
\end{equation}
This boundary is pseudoconvex, has bracket finite type $2m$ and \eqref{1.4} is satisfied. Thus Theorem~\ref{t1.1} applies and we have \eqref{1.5}. This is more general than  \cite{BDKT06} where it is assumed \eqref{b}.
Thus, for example, for the domain graphed by $g$ with
$$
g=|z|^{2{(m-1)}}x^2h(z)\qquad\T{for $h>0$ and $h_{1\bar1}>0$},
$$
we have \eqref{a} though the second of \eqref{b} is never true, not even for $h\equiv1$. For general $h$, neither of \eqref{b} is fulfilled.
\ee
There is a result for sum of squares which stays  close to Theorem~\ref{p1.1} and in particular only assumes finite type without the additional hypothesis \eqref{1.4}. This requires to modify the Kohn-Laplacian as 
$$
\tilde\Box^k=\Lambda_{\di_t}^{-2l}L\bar L+L|z|^{2k}\bar L+\bar L|z|^{2k}L,
$$
where $\Lambda^{-2l}_{\di_t}$ is the standard pseudodifferential operator of order $-2l$ in $t$.
\bt
\Label{t1.2}
Let $\{L,\,\bar L\}$ have type $2m$; then
\begin{equation}
\Label{extranew}
\no{\zeta u}^2_s\simleq \no{\zeta'\tilde\Box^ku}^2_{s+2l}+\no{u}^2_0.
\end{equation}
\et
\vskip0.2cm
Some references to current literature are in order. Hypoellipticity in presence of infinite degeneracy has been intensively discussed in recent years. The ultimate level to which the problem is ruled by estimates, are superlogarithmic estimates (Kusuoka and Strooke \cite{KS85}, Morimoto \cite{M87} and Kohn \cite{K02}). Related work is also by Bell and Mohammed \cite{BM95} and Christ \cite{Ch02}. Beyond the level of estimates are the results by Kohn \cite{K00} which develop, in a geometric framework, an early result by Fedi \cite{F71}: the point here is that the degeneracy is confined to a real curve transversal  to the system. This explains also why if the set of degeneracy is big, superlogarithmicity becomes in certain cases necessary (\cite{M87} and \cite{Ch02}). In all these results, however, there is somewhat a gain of derivatives (such as sublogarithmic). The simplest example of hypoellipticity without gain (nor loss) is $\Box_b+\lambda\,\T{id}$, $\lambda>0$ where $\Box_b$ is the Kohn-Laplacian of $\Re w=|z|^2$ (cf. Stein \cite{St82} where the bigger issue of the analytic-hypoellipticity is also addressed). 
As for  loss of derivatives, the phenomenon has been discovered by Kohn in \cite{K05} and further developped by Bove, Derridj, Kohn and Tartakoff in \cite{BDKT06}. Additional contribution is, among others, by Parenti and Parmeggiani \cite{PP05} and Tartakoff \cite{T06}.

\section{Sums of squares in $\R^{2n+1}$ for $n>1$}
\Label{s2}
We restate in higher dimension the results of Section~\ref{s1}; we can better appreciate the different role which is played by the finite type with respect to  \eqref{1.4}. The containt of this section is a direct consequence of the results of Section~\ref{s1} (plus ellipticity and maximal hypoellipticity related to microlocalization) and therefore it does not need a specific proof.
In $\C^n\times \R_t$ we start, as in Section~\ref{s1}, from $L_1=\di_{z_1}-ig_1(z_1)\di_t$  and  complete $L_1$ to a system of smooth complex vector fields in a neighborhood of $0$
$$
L_j=\partial_{z_j}-ig_j(z)\partial_t,\,\,j=1,...,n\qquad \T{for $g_j|_0=0$}.
$$
For a system of vector fields, we denote by ${\mathcal Lie}_{2m}$ the span of commutators of order $\leq 2m-1$ belonging to the system.
We have
 $\no{u^0}^2_1\simleq \sumj\no{\bar L_ju^0}^2_0+\NO{u}_0$ and, if for some index $j$, say $j=1$, $\di_t\in{\mathcal Lie}_{2m_1}\{ L_1,\,\bar L_1\}$, then $\no{u^-}^2_{\frac1{2m_1}}\simleq \sumj\no{\bar L_ju}^2_0+\NO{u}_0$ (cf. the end of Section~\ref{s3}). Summarizing up, if we only have \eqref{1.2} for $u^+$, we get, for the full $u$ and with $l$ replaced by $l_1=\frac{k_1}{2m_1}$:
\begin{equation}
\Label{1.2ter}
\no{\zeta u}^2_{s}\simleq \left(\no{\zeta'\bar L_1 u}^2_{s-\frac1{2m_1}}+\no{\zeta'z_1^{k_1}\bar L_1 u}^2_{s+l_1}
+\no{\zeta' z_1^{k_1}L_1 u}^2_{s+l_1}\right)+\underset{j=2}{\overset n\sum}\no{\bar L_j u}^2_{s-\frac1{2m_1}}+\no{u}^2_0.
\end{equation}
We assume that each coefficient satisfy $g_j=\partial_{z_j}g$ for a real function $g=g(z),\,\,z=(z_1,...,z_n)\in\C^n$ 
and denote by $\Bbb L$ the bundle spanned by the $L_j$'s. We note that this defines a CR structure because, 
on account of $g_{i\bar j}=g_{j\bar i}$,
$$
\T{ $\Bbb L$ is involutive}.
$$
Also, this structure is of hypersurface type in the sense that
$$
T(\C^n_z\times \R_t)=\Bbb L\oplus \overline{\Bbb L}\oplus \R\partial_t.
$$
Note that, in fact, the $L_j$'s commute; therefore, the Levi form is defined directly by $[L_i,\bar L_j]=g_{ij}\partial_t$, without passing to the quotient modulo $\Bbb L\oplus \overline{\Bbb L}$.
We also assume that the Levi form $(g_{i\bar j})$ is positive semidefinite; in particular $g_{j\bar j}\geq 0$ for any $j$. (Geometrically, this means that the hypersurface $\Im w=g$ graphed by $g$, is pseudoconvex.) 
We choose $\kappa=(k_1,...,k_n)$ and define the perturbed Kohn-Laplacian
$$
\Box^\kappa=\sumj L_j\bar L_j+\bar L_j|z_j|^{2k_j}L_j.
$$ 
\bt
\Label{t0.1}
Assume that for any $j$, $\di_t\in{\mathcal Lie}_{2m_j}\{ L_j,\,\bar L_j\}$, and that
\begin{equation}
\Label{0.3}
|g_j|\simleq |z_j|g_{j\bar j}\,\,\,\T{and}\,\,\, |g_{j\bar j\bar j}|\simleq |z_j|^{-1}g_{j\bar j}\quad \T{for any $j=1,...,n$.}
\end{equation}
Define $l_j:=\frac {k_j-1}{2m_j}$ and put $l=\underset j \max \frac {k_j-1}{2m_j}$; then
\begin{equation}
\Label{1.4new}
\no{\zeta u}^2_s\simleq \no{\zeta' \Box^\kappa u}^2_{s+2l}+\no{u}^2_0.
\end{equation}
\et
The proof of Theorem~\ref{t0.1} and Theorem~\ref{t0.2} below, are just a variation of those of the twin Theorems \ref{t1.1} and \ref{t1.2}.
We define now
$$
\tilde \Box^\kappa=\sumj \left(\Lambda_{\di_t}^{-2l_j}L_j\bar L_j+\sumj L_j|z_j|^{2k_j}\bar L_j+\bar L_j|z_j|^{2k_j}L_j\right).
$$
\bt
\Label{t0.2}
Assume that for any $j$, $\di_t\in{\mathcal Lie}_{2m_j}\{ L_j,\,\bar L_j\}$; then
\begin{equation}
\Label{0.5}
\no{\zeta u}^2_s\simleq \no{\zeta'\tilde\Box^ku}^2_{s+2l}+\no{ u}^2_0.
\end{equation}
\et

\section{Preliminaries}
\Label{s3}
We identify $\C\times\R$ to $\R^3$ with coordinates $(z,\bar z,t)$ or $(\Re z,\Im z, t)$.
We denote by $\xi=(\xi_z,\xi_{\bar z},\xi_t)$ the variables dual to $(z,\bar z,t)$, by $\Lambda^s_\xi$ the standard symbol $(1+|\xi|^2)^{\frac s2}$, and by 
$\Lambda^s_\partial$ the pseudodifferential operator with symbol $\Lambda_\xi^s$; this is defined by $\Lambda^s_\partial(u)=\mathcal F^{-1}(\Lambda_\xi
^s\mathcal F(u))$ where $\mathcal F$ is the Fourier transform. 
We consider the full (resp. totally real) $s$-Sobolev norm $\no{u}_s:=\no{\Lambda^s_\partial u}_0$ (resp. $\no{u}_{\R, \,s}:=\no{\Lambda_{\partial_t}^su}_0
$).
In $\R^{3}_\xi$, we consider a conical partition of the unity $1=\psi^++\psi^++\psi^0$ where $\psi^\pm$ have support in a neighborhood of the axes $\pm\xi_t$ and $\psi^0$ in a neighborhood of the plane $\xi_t=0$, and introduce a decomposition of the identity $\T{id}=\Psi^++\Psi^-+\Psi^0$ by means of $\Psi^{\overset\pm0}$, the pseudodifferential operators with symbols $\psi^{\overset\pm0}$; we accordingly write $u=u^++u^-+u^0$.
Since $|\xi_z|+|\xi_{\bar z}|\simleq \xi_t$ over $\T{supp}\,\psi^+$, then  $\no{u^+}_{\R, \,s}=\no{u^+}_s$.

We carry on the discussion by describing the properties of commutation of the vector fields  $L$ and $\bar L$ for $L=\di_z-ig_1(z)\di_t$.
The crucial equality is 
\begin{equation}
\Label{2.1}
\no{L  u}^2=([L ,\bar L ]u,u)+\no{\bar L  u}^2,\quad u\in C^\infty_c,
\end{equation}
which is readily verified by integration by parts. 
Note here that errors coming from derivatives of coefficients do not occur since $g_1$ does not depend on $t$. 
Recall that $[L,\bar L]=g_{1\bar1}\di_t$; this implies
\begin{equation}
\Label{supernova}
|(g_{1\bar 1}\di_t u,u)|\simleq s.c. \no{\di_tu}^2+l.c.\no{u}^2.
\end{equation}
We have
\begin{equation}
\begin{split}
\Label{2.2}
\NO{u^0}_1&\simleq \NO{\bar Lu^0}+\NO{Lu^0}+\NO{u}
\\
&\leq 2\NO{\bar Lu^0}+ sc \NO{\di_tu^0}+lc\NO{u}.
\end{split}
\end{equation}
To check \eqref{2.2},  we point our attention to the estimate for operator's symbols $(1+|\xi|^2)|\alpha|^2\simleq |\alpha|^2+|\sigma(\bar L)\alpha|^2+|\sigma(L)\alpha|^2$ ($\alpha$ complex) over $U\times \T{supp}\,\psi^0$ for a neighborhood $U$ of $0$; in addition to the fact that $[L,\Psi^0]$ is of order $0$, this yields the first inequality of \eqref{2.2}. The second follows from \eqref{2.1} combined with \eqref{supernova}. As  for $u^-$, since $g_{11}\sigma(\di_t)<0$ over $\T{supp}\,\psi^-$, then
\begin{equation*}
(g_{11}\di_tu^-,u^-)=-\left|(g_{11}\Lambda_{\di_t}u^-,u^-)\right|.
\end{equation*}
Thus \eqref{2.1} implies $\no{L u^-}\leq \no{\bar L  u^-}$ (the second inequality in \eqref{2.3} below). Suppose now that $\{L,\,\bar L\}$ have type $2m$; this yields the first inequality below which, combined with the former, yields
\begin{equation}
\Label{2.3}
\begin{split}
\no{u^-}^2_{\frac1{2m}}&\simleq \no{L u^-}^2_0+\no{\bar Lu^-}^2_0+\NO{u}_0
\\
&\simleq \no{\bar L u^-}^2_0+\NO{u}_0.
\end{split}
\end{equation}
In conclusion, only estimating $u^+$ is relevant. For this purpose, we have a useful statement
\bl
\Label{l2.2}
Let $|[L ,\bar L ]|^{\frac12}$ be the operator with symbol $|g_{11}|^{\frac12}\Lambda^{\frac12}_{\xi_t}$; then 
\begin{equation}
\Label{2.5}
\no{|[L , \bar L ]|^{\frac12}u^+}^2\leq \no{L u^+}^2+\no{\bar L  u^+}^2.
\end{equation}
\el
\bpf
From \eqref{2.1} we  get
\begin{equation*}
|([L ,\bar L ]u,u)|\leq \no{L u}^2+\no{\bar L u}^2.
\end{equation*}
The conclusion then follows from
$$
[L ,\bar L ]=|[L ,\bar L]|\qquad\T{over $\T{supp}\,\psi^+$}.
$$

\epf

We pass to a result about intepolation which plays a central role in our discussion.
\bl
\Label{l2.1}
Let $f=f(z)$ be smooth and satisfy $f(0)=0$. Then for any $\rho$, $r$, $n_1$ and $n_2$ with $0<n_1\leq r$, $n_2>0$
\begin{equation}
\Label{2.4}
\no{f^ru}^2_0\simleq sc \no{f^{r-n_1}u}^2_{\R,\,-n_1\rho}+ lc \no{f^{r+n_2}u}^2_{\R,\,n_2\rho}.
\end{equation}
\el
\bpf
Set  $A:=\Lambda^\rho_{\di_t}f$; interpolation for the pseudodifferential operator $A$ yields
\begin{equation*}
\begin{split}
\no{f^ru}^20&=\no{(\Lambda^\rho f)^ru}^2_{\R\,-\rho r}
\\
&=(\Lambda^{\rho(r-n_1)}f^{r-n_1},\Lambda^{\rho(r+n_1)}f^{r+n_1})_{-\rho r}
\\
&=(\Lambda^{-\rho n_1}f^{r-n_1},\Lambda^{\rho n_1}f^{r+n_1})_0\simleq sc\no{f^{r-n_1}}^2_{\R,\,-n_1\rho}+ lc \no{f^{r+n_1}u}^2_{\R,\,n_1\rho}.
\end{split}
\end{equation*}
This proves the lemma for $n_2=n_1$; the general conclusion is obtained by iteration.

\epf
We have now a result about factors in a scalar product.
\bl
\Label{l2.3}
Let $h=h(z)$ satisfy $|h|\leq |h_1||h_2|$ and take $f=f(z,t)$ and $g=g(z,t)$. Then
\begin{equation}
\Label{2.5,5}
|(f,hg)|_{\R,\,s}\simleq \no{fh_1}^2_{\R,\,s}+\no{gh_2}^2_{\R,\,s}.
\end{equation}
\el
\bpf
We use the notation ${\mathcal F}_t $ for the partial Fourier transform with respect to $t$ and $d\lambda$ for the element of volume in $\C_z\simeq \R^{2}_{\Re z,\Im z}$. The lemma follows from the following sequence of inequalities in which the crucial fact is that $h,\,h_1$ and $h_2$ are constant in the integration in $\xi_t$:
\begin{equation*}
\begin{split}
|(f,hg)_{\R,\,s}|&=\left|\int_{\R^{2}}\left(\int_{\R^1_{\xi_t}}\Lambda^{2s}_{\xi_t}{\mathcal F}_t (f)h{\mathcal F}_t (g)d\xi_t\right)d\lambda\right|
\\
&\leq \int_{\R^{2}}\left(\int_{\R^1_{\xi_t}}\Lambda^{2s}_{\xi_t}|{\mathcal F}_t (f)h_1h_2{\mathcal F}_t (g)|d\xi_t\right)d\lambda
\\
&\leq \int_{\R^{2}}\left(\int_{\R^1_{\xi_t}}\Lambda^{2s}_{\xi_t}|{\mathcal F}_t (f)h_1|^2d\xi_t\right)d\lambda+\int_{\R^{2}}\left(\int_{\R^1_{\xi_t}}\Lambda^{2s}_{\xi_t}|{\mathcal F}_t (g)h_2|^2d\xi_t\right)d\lambda
\\
&\underset{\T{Plancherel}}=\no{fh_1}^2_{\R,\,s}+\no{gh_2}^2_{\R,\,s}.
\end{split}
\end{equation*}

\epf

We say a few words for the case of higher dimension. In $\C^{n}_{z_1,...,z_n}\times \R_t$, we consider a full system $L_j=\di_{z_j}-ig_j\di_t,\,\,j=1,...,n$ with $g_j|_0=0$. The same argument used in proving \eqref{2.2} yields
\begin{equation}
\Label{2.2bis}
\NO{u^0}_1\simleq \sumj\NO{\bar L_ju^0}+\NO{u}.
\end{equation}
Similarly as above, we have $\no{L_j u^-}^2\leq \no{\bar L_j  u^-}^2+\NO{u}$ for any $j$. Then, if at least one index $j$, say $j=1$, the pair $\{L_1,\,\bar L_1\}$ has type $m=m_1$, we get, in the same way as in \eqref{2.3}
$$
\no{u^-}^2_{\frac1{2m}}\simleq \sumj\no{\bar L_j u^-}^2+\NO{u}.
$$
Again, only estimating $u^+$ is therefore relevant.

\section{Proof of Theorem~\ref{p1.1}, Theorem~\ref{t1.1} and Theorem~\ref{t1.2}}
\Label{s4}
\noindent
{\it Terminology.}\hskip0.2cm In an estimate we call ``good" a term in the right side (upper bound). 
We call ``absorbable"  a term that we encounter in the course of the estimate and which comes as a fraction (small constant or sc) of a former term.
If  cut-off  are involved in the estimate, and in the right side the cut-off can be expanded, say passing from $\zeta$ to $\zeta'$, we call ``neglectable" a term which comes with lower Sobolev index and possibly with a bigger cut-off. Neglectable is meant with respect to the initial (left-hand side) term of the estimate, to further terms that one encounters and even to extra terms provided that they can be estimated by ``good". These latter are sometimes artificially added to expand the range of ``neglectability".

\vskip0.2cm
\noindent
{\it Proof of Theorem~\ref{p1.1}.} \hskip0.2cm According to \eqref{2.2} and \eqref{2.3}, it suffices to prove \eqref{1.2} for $u=u^+$; so, throughout the proof we write $u$ but mean $u^+$. Also, we use the equivalence, over $u^+$, between the totally real $\no{\cdot}_{\R,\,s}-$ with the full 
$\no{\cdot}_s$-Sobolev norm; the specification of the norm will be omitted.  Moreover, we  can use a cut-off $\zeta=\zeta(t)$ in $t$ only. In fact, for a cut-off $\zeta=\zeta(z)$ we have $[L ,\zeta(z)]=\dot\zeta$ and $\dot\zeta\equiv0$ at $z=0$. On the other hand, $z^{k}L \sim L $ outside $z=0$ which yields  \eqref{3.1} below (so that we have gain, instead of loss).   Recall in fact that we are assuming that $M$ has type $2m$. It is classical that the tangential vector fields $L$ and $\bar L$ satisfy $\frac1{2m}$-subelliptic estimates, that is, the first inequality in the estimate below. In combination with \eqref{2.1} which implies the second inequality below, we get
\begin{equation}
\Label{3.1}
\begin{split}
\no{\zeta u}^2_{s}&\simleq \no{\zeta \bar L u}^2_{s-\frac1{2m}}+\no{\zeta L u}^2_{s-\frac1{2m}}+\no{\zeta'u}^2_{s-\frac1{2m}}
\\
&\simleq \no{\zeta \bar L u}^2_{s-\frac1{2m}}+\no{\zeta \left|[L ,\bar L ]\right|^{\frac12}u}^2_{s-\frac1{2m}}+\no{\zeta' u}^2_{s-\frac1{2m}}.
\end{split}
\end{equation}
Remark that $\NO{\zeta'u}_{s-\frac1m}$ (for a new $\zeta'$) takes care of the error $\NO{\zeta'\bar Lu}_{s-\frac1{2m}-1}$ coming from $[\Lambda^{2s-\frac1m},\zeta']$. 
Now, remember that $[L ,\bar L ]=g_{1\bar1}\di_t$ without error terms, that is, combinations of $L$ and $\bar L$; recall also that $g_{1\bar1}\geq0$. We get
\begin{equation}
\Label{3.2}
\begin{split}
\no{\zeta \left|[L ,\bar L ]\right|^{\frac12}u}^2_{s-\frac1{2m}}&\sim\no{\zeta g_{1\bar1}^{\frac12}\Lambda^{\frac12}_{\partial_t}u}^2_{s-\frac1{2m}}
\\
&\simleq sc\, \no{\zeta u}^2_{s}+ lc \,\no{\zeta g_{1\bar1}^{\frac12+\frac {k}{2(m-1)}}\Lambda_{\di_t}^{\frac12}u}^2_{s+l}
\\
&\simleq \T{absorbable}+\no{\zeta g_{1\bar1}^{\frac12}z^{k}\Lambda^{\frac12}_{\di_t}u}^2_{s+l}
\\
&=\T{absorbable}+\no{\zeta \left|[L ,\bar L ]\right|^{\frac12}z^ku}^2_{s+l}
\\
&\leq\T{absorbable}+\no{\zeta L (z^{k}u)}^2_{s+l}+\no{\zeta \bar L (z^{k}u)}^2_{s+l}+\no{\zeta'z^{k}u}^2_{s+l},
\end{split}
\end{equation}
where the first ``$\sim$" is a way of rewriting the commutator, the second ``$\simleq$" follows from Lemma~\ref{l2.1}
(under the choice $n_1=m-1$, $n_2=k$, $r=m-1$, $\rho=\frac1{2m}$ and $f=g_{1\bar1}^{\frac1{2(m-1)}}$),
the third ``$\simleq$" follows from $|g_{1\bar1}|\simleq |z|^{2(m-1)}$, the fourth ``$=$" is obvious and the last ``$\simleq$" follows from Lemma~\ref{l2.2}.  We go now to estimate, in the last line of \eqref{3.2}, the two terms $\no{\zeta \bar L (z^{k}u)}^2_{s+l}$ and $\no{\zeta 'z^{k}u}^2_{s+l}$. We start from
\begin{equation}
\Label{3.2,5}
\no{\zeta  L (z^{k}u)}^2_{s+l}\leq \no{\zeta z^{k} L u}^2_{s+l}+\no{\zeta z^{k-1}u}^2_{s+l},
\end{equation}
where the last term is produced by the commutator $[L ,z^{k}]$. By writing, in the scalar product, once $z^{k-1}$ and once $[L ,z^{k}]$, we get
\begin{equation}
\Label{3.3}
\begin{split}
\no{\zeta z^{k-1}u}^2_{s+l}&=(\zeta z^{k-1}u,\zeta [L ,z^{k}]u)_{s+l}
\\
&=(\zeta z^{k-1}u,\zeta z^{k}L u)_{s+l}+(\zeta z^{k-1}u,\zeta L  z^{k}u)_{s+l}.
\end{split}
\end{equation}
Now,
\begin{equation}
\Label{3.4}
\begin{cases}
(\zeta z^{k-1}u,\zeta z^{k}L u)_{s+l}\leq \underset{\T{absorbable}}{\underbrace{sc \no{\zeta z^{k-1}u}^2_{s+l}}}+\underset{\T{good}}{\underbrace{\no{\zeta z^{k}L u}^2_{s+l}}}
\\
\begin{split}
(\zeta z^{k-1}u,\zeta L z^{k}u)_{s+l}&=(\zeta z^{k-1}\bar L u,\zeta z^{k}u)_{s+l}+(\zeta z^{k-1},\zeta'z^{k}u)_{s+l}
\\
&\underset{\T{good}}{\underbrace{\simleq\no{\zeta z^{k}\bar L  u}^2_{s+l}}}+\underset{\T{absorbable}}{\underbrace{sc \no{\zeta z^{k-1}u}^2_{s+l}}}+\no{\zeta'z^{k}u}^2_{s+l}.
\end{split}
\end{cases}
\end{equation}
Thus $\no{\zeta z^{k-1}u}^2_{s+l}$ has been estimated by $\no{\zeta'z^{k}u}^2_{s+l}$. What we have obtained so far is
\begin{equation}
\Label{estimate}
\no{\zeta u}^2_s\simleq \no{\zeta \bar L u}^2_{s }+\no{\zeta z^k\bar L u}^2_{s+l }+\no{\zeta z ^{k}L u}^2_{s+l }+\no{\zeta'z^ku}^2_{s+l}+\no{\zeta'u}^2_{s-\frac1{2m}}.
\end{equation}
Note that in this estimate, the terms coming with $L$ and $\bar L$ carry the same cut-off $\zeta$ as the left side; it is in this form that Theorem~\ref{p1.1} will be applied for the proof of Theorems~\ref{t1.1} and \ref{t1.2}.
Instead, to conclude the proof of Theorem~\ref{p1.1}, we have to go further with the estimation of  $\no{\zeta'z^{k}u}^2_{s+l}$ (which also provides the estimate of the last term in \eqref{3.2}). We have, by subelliptic estimates
\begin{equation}
\Label{3.5}
\no{\zeta'z^{k}u}^2_{s+l}\simleq \no{\zeta'L z^{k}u}^2_{{s+l}-\frac1{2m}}+\no{\zeta'\bar L  z^{k}u}^2_{{s+l}-\frac1{2m}}+\no{\zeta''z^{k}u}^2_{{s+l}-\frac1{2m}}.
\end{equation}
To $\no{\zeta'Lz^{k}u}^2_{{s+l}-\frac1{2m}}$ we apply \eqref{3.2,5} with ${s+l}$ replaced by ${s+l}-\frac1{2m}$. In turn, $\no{\zeta' z^{k-1}u}^2_{s+l-\frac1{2m}}$ can be estimated, by \eqref{3.3}, \eqref{3.4} and \eqref{3.5} with Sobolev indices all lowered from $s+l$ to $s+l-\frac1{2m}$, by means of 
``good" + ``absorbable" + $\no{\zeta''z^ku}^2_{s+l-\frac1{2m}}$. (In fact, ``good" even comes with lower index.)
 The conclusion \eqref{1.2} follows from induction over $j$ such that $\frac j{2m}\geq s+l$.
This completes the proof of Theorem~\ref{p1.1}.

\hskip15cm$\Box$

\vskip0.3cm
\noindent
{\it Proof of Theorem~\ref{t1.2}.}
\hskip0.2cm We first prove Theorem~\ref{t1.2} instead of Theorem~\ref{t1.1} because it is by far easier. As it has already been remarked in Section 1, it suffices to prove the theorem for $u=u^+$. Also, in this case, the full norm can be replaced by the totally real norm.
So we write $u$ for $u^+$ and $\no{\cdot}_s$ for $\no{\cdot}_{\R,\,s}$; however, in some crucial passage where Lemma~\ref{l2.3} is on use, it is necessary to point attention to the kind of he norm.
We start from \eqref{estimate}; note that, for this estimate to hold, only finite type is required.
We begin by noticing that the last term of \eqref{estimate} is neglectable.
We then rewrite the third term in the right of \eqref{estimate} as
\begin{equation}
\Label{3.8}
(\zeta z^{k}L u,\zeta z ^{k}L u)_{s+l }=(\zeta\bar L |z |^{2k}L u,\zeta u)_{s+l }+(\zeta z ^{k}L u,\zeta'z ^{k}u)_{s+l },
\end{equation}
where we recall that we are using the notation $l =\frac {k-1}{2m}$.
(Note that the commutator $[L ,\zeta]$ is not just $\zeta'$ but comes with an additional factor $g_1$, the coefficient of $L $; but we disregard this contribution here though it will play a crucial role in the proof of Theorem~\ref{t1.1}.)
We keep the first term in the right of \eqref{3.8} as it stands and put together with the similar term coming from the first term in the right of 
\eqref{estimate} to form $\tilde\Box^\kappa$. 
We then apply Cauchy-Schwartz inequality and estimate the first term by $\no{\zeta\tilde\Box^\kappa u}^2_{s+2l }+sc\no{\zeta u}^2_s$.
As for the second term in the right of \eqref{3.8}, it can be estimated, via Cauchy-Schwartz, by  $sc \no{\zeta z ^{k}Lu}^2_{s+l }+lc\no{\zeta'z^ku}_{s+l}$. To this latter, we apply subelliptic estimates 
\begin{equation}
\Label{3.9}
\no{\zeta'z ^{k}u}^2_{s+l }\simleq \no{\zeta'z^k\bar L u}^2_{s+l-\frac1{m}}+\no{\zeta'z ^{k}L u}^2_{s+l -\frac1{2m}}+\no{\zeta'z ^{k-1}u}^2_{s+l -\frac1{2m}}+\no{\zeta'' z ^{k}u}^2_{s+l -\frac1{2m}}.
\end{equation}
For the third term in the right, recalling  \eqref{3.3} and \eqref{3.4}, we get
\begin{equation}
\Label{3.10}
\no{\zeta' z ^{k-1}u}^2_{s+l -\frac1{2m}}\simleq\T{neglectable}+\no{\zeta''z ^{k}u}^2_{s+l -\frac1{2m}}.
\end{equation}
Thus $\no{\zeta'z ^{k}u}^2_{s+l }$ is controlled by induction over $j$ with $\frac j{2m}\geq s+l $. (Recall, once more, that ``good" is stable under passing from $\zeta'$ to $\zeta''$.)
We notice that combination of \eqref{3.9} and \eqref{3.10} shows that  $\no{\zeta'z ^{k}u}^2_{s+l }$ is neglectable.
We pass to $\no{\zeta'z ^{k}\bar L u}^2_{s+l}$, the second term in the right of \eqref{estimate} and observe that it  can be treated exactly in the same way as the third (with $L$ instead of $\bar L$).
We end with the first which does not carry the loss $l$; we have
\begin{equation}
\Label{3.11}
\begin{split}
\no{\zeta\bar L  u}^2_{s}&=(\zeta L \bar L u,\zeta u)_{s }+(\zeta \bar L u,\zeta'g_1u)_{s }
\\
&=(\Lambda^{2l}\Lambda^{-2l}L\bar L u,\zeta u)_s+(\zeta \bar L u,\zeta'g_1u)_{s }.
\end{split}
\end{equation}
The first term in the right combines to form $\tilde \Box^k$. As for the second,
we notice that $|g_1|\simleq |z|$ and therefore applying Lemma~\ref{l2.1} for $n_1=k-1$ and $n_2=1$
$$
(\zeta \bar L u,\zeta'g_1u)_{s }\leq sc \no{\zeta \bar Lu}^2_s+ lc (\no{\zeta'z^ku}^2_{s+l}+\no{\zeta' u}^2_{s-\frac1{2m}}).
$$
The first term in the right is absorbable, the last neglectable, the midle has already been proved to be neglectable by subelliptic estimates \eqref{3.9}.
This completes the proof.

\hskip15cm$\Box$

\vskip0.2cm
\noindent
{\it Proof of Theorem~\ref{t1.1}.} \hskip0.2cm 
As before, we prove the theorem for $u=u^+$ and write $\no{\cdot}_s$ for $\no{\cdot}_{\R,\,s}$ though, in some crucial passage, it is necessary to point the attention to the kind of the norm.
Raising Sobolev indices, we
 rewrite \eqref{estimate} in a more symmetric fashion as
\begin{equation}
\Label{3.12}
\no{\zeta u}^2_s\simleq \no{\zeta\bar L u}^2_{s+l }+\no{\zeta z ^{k}L u}^2_{s+l }+\no{\zeta'u}^2_{s-\frac1{2m}}+\no{\zeta'z^{k}u}^2_{s+l }.
\end{equation}
We handle all terms in the right as in Theorem~\ref{t1.1} except from the first which comes now with the loss $s+l$.
We point out that to control these terms,  only  finite type has been used. Instead, to control the remaining term, we need the additional hypothesis \eqref{1.4}. We have
\begin{equation}
\Label{3.13}
\no{\zeta\bar L  u}^2_{s+l }=(\zeta L \bar L u,\zeta u)_{s+l }+(\zeta \bar L u,\zeta'g_1u)_{s+l }. 
\end{equation}
The first term combines to form $\Box^k$. As for the second, we recall the estimate $|g_1|\simleq |z|g_{1\bar1}$ and apply Lemma~\ref{l2.3} for $h=zg_{1\bar1},\, h_1=g_{1\bar1}^{\frac12}$ and $h_2=zg_{1\bar1}^{\frac12}$ to get
\begin{equation}
\Label{3.14}
|(\zeta\bar L u,\zeta'g_1u)|_{s+l }\leq sc \no{\zeta g_{1\bar1}^{\frac12}\bar L u}^2_{s+2l }+ lc \no{\zeta'zg_{1\bar1}^{\frac12}u}^2_{s }.
\end{equation}
In the estimate above, we point our attention to the fact that  the norms that we are considering are totally real norms (though we do not keep track in our notation) and therefore Lemma~\ref{l2.3} can be applied. 
We start by estimating the second term in the right. By Lemma \ref{l2.2} and next, Lemma~\ref{l2.1} for $n_1=1$, $n_2=k-1$
\begin{equation}
\Label{3.15}
\begin{split}
\no{\zeta' g_{1\bar1}^{\frac12}	zu}^2_s&\simleq \no{\zeta'z Lu}^2_{s-\frac12}+\no{\zeta' z\bar Lu}^2_{s-\frac12}+\T{neglectable}
\\
&\leq  \no{z^k\zeta'Lu}^2_{s-\frac12+l}+  \no{\zeta' Lu}^2_{s-\frac12-\frac1{2m}}+\no{\zeta' z\bar Lu}^2_{s-\frac12}+\T{neglectable},
\end{split}
\end{equation}
where $\T{neglectable}$ comes from the commutators $[L,z]$ and $[L,\zeta']$. Also, the first  term in the second line of \eqref{3.15} is neglectable.
As for the second term, we have, by \eqref{2.1}
\begin{equation}
\Label{3.16}
\no{\zeta' Lu}^2_{s-\frac12-\frac1{2m}}\simleq \no{\zeta' g_{1\bar1}^{\frac12}u}^2_{s-\frac1{2m}}+\no{\zeta'\bar Lu}^2_{s-\frac12-\frac1{2m}}+\T{neglectable}.
\end{equation}
Since both  terms in the right of \eqref{3.16} are neglectable, we conclude that $\no{\zeta'zg_{1\bar1}^{\frac12}u}^2_s$ itself is neglectable. 
From now on, we follow closely the track of \cite{BDKT06}.
We pass to consider the last and most difficult term to estimate, that is, the first in the right of \eqref{3.14}. Along with this term, that we denote by (a), we introduce three additional terms; we set therefore
\begin{equation*}
\begin{cases}
(a):=\no{\zeta g_{1\bar1}\bar L u}^2_{s+2l },& (b):=\no{\zeta z^{2k-1}g_{1\bar1}^{\frac12}u }_{s+2l},
\\
(c):=\no{\zeta z^{2k-1}Lu}^2_{s+2l-\frac12},& (d):=\no{L\zeta\bar Lu}^2_{s+2l-\frac12}.
\end{cases}
\end{equation*}
Because of these additional terms, that we are able to estimate, ``neglectable" and ``absorbable" take an extended range. We first show that (b) is controlled by (c). This is apparently as in \cite{BDKT06} first half of 5.3 but more complicated because our (b) and (c) are different from their $(LHS)_5$ and $(LHS)_6$ respectively. Now, by Lemma~\ref{l2.2} we get
\begin{equation*}
(b)\simleq (c)+\no{\zeta'z^{2k-1}g_1u}^2_{s+2l-\frac12}+\no{\zeta z^{2k-2}u}^2_{s+2l-\frac12}+\T{neglectable},
\end{equation*}
where the central terms in the right come from $[L,\zeta]$ and $[L,z^{2k-1}]$ respectively, and where neglectable, with respect to (a), is the term which involves $\bar Lu$ and which comes lowered by $-\frac12$. The first of the central terms is neglectable with respect to (b). As for the second, we have, using the notation $\#=s+2l-\frac12-\frac1{2m}$
\begin{equation*}
\Label{3.17}
\begin{split}
\underset{(i)}{\underbrace{\no{\zeta z^{2k-2}u}^2_{s+2l-\frac12}}}&\simleq \underset{(ii)}{\underbrace{\no{\zeta z^{2k-2}Lu}_\#}}+\underset{(iii)}{\underbrace{\no{\zeta z^{2k-2}\bar Lu}^2_{\#}}}
\\
&+\underset{(iv)}{\underbrace{\no{\zeta' z^{2k-2}g_1u}^2_\#}}+\underset{(v)}{\underbrace{\no{\zeta z^{2k-3}u}^2_\#}},
\end{split}
\end{equation*}
where the two terms of the second line come from $[L,\zeta]$ and $[L,z^{2k-2}]$ respectively. First, (iv) is neglectable with respect to (i). Next, using Lemma~\ref{l2.1} for $n_1=2k-2$ and $n_2=1$
\begin{equation*}
(ii)\simleq \underset{(ii)_1}{\underbrace{sc \no{\zeta Lu}^2_{\#-\frac{2k-2}{2m}}}}+\underset{(ii)_2}{\underbrace{lc \no{\zeta z^{2k-1}Lu}^2_{\#+\frac1{2m}}}}.
\end{equation*}
Note that $\#-\frac{2k-2}{2m}=s-\frac12-\frac1{2m}$ and $\#+\frac1{2m}=s+2l-\frac12$; thus
 $(ii)_1$ is absorbed by \eqref{3.16} and $(ii)_2$ is estimated by (c). Next, by Lemma~\ref{l2.1} for $n_1=2k-3$ and $n_2=1$
 $$
 (v)\simleq lc\underset{(v)_1}{\underbrace{\no{\zeta u}^2_{\#-\frac{2k-3}{2m}}}}+ sc \underset{(v)_2}{\underbrace{\no{\zeta z^{2k-2}u}^2_{\#+\frac1{2m}}}}.
$$
We have $\#-\frac{2k-3}{2m}=s-\frac12$ and, again, $\#+\frac1{2m}=s+2l-\frac12$; thus
  $(v)_1$ is neglectable with respect to $\no{\zeta u}^2_s$, the term in the left of the estimate, and $(v)_2$ is absorbed by (i). Finally, by \eqref{2.1} 
  \begin{equation*}
(iii)\leq \underset{(iii)_1}{\underbrace{\no{\zeta z^{2k-2}Lu}_{\#}}}+\underset{(iii)_2}{\underbrace{\no{\zeta z^{2k-2}g_{1\bar1}^{\frac12}u}^2_{\#+\frac12}}}.
\end{equation*}
Now, applying Lemma~\ref{l2.1} for $n_1=k-2$, $n_2=1$ in the first line below and $n_1=2k-2$ and $n_2=1$ in the second respectively, we get
\begin{equation*}
\begin{cases}
(iii)_1\simleq \underset{(c)}{\underbrace{\no{\zeta z^{2k-1}Lu}_{s+2l-\frac12}}}+\underset{\T{neglectable w.r.to $\no{\zeta z^k Lu}^2_{s+l}$}}{\underbrace{\no{\zeta z^kLu}_{s+l-\frac12}}}
\\
(iii)_2\simleq sc \underset{(b)}{\underbrace{\no{\zeta z^{2k-1}g_{1\bar1}^{\frac12}u}^2_{s+2l}}}+lc \underset{\T{neglectable w.r.to $\no{\zeta u}_s$}}{\underbrace {\no{\zeta g_{1\bar1}^{\frac12}u}^2_{s-\frac1{2m}}}}.
\end{cases}
\end{equation*}
Summarizing up,
$$
(b)\simleq (c)+\T{neglectable}.
$$
We have to show now that 
\begin{equation*}
\begin{cases}
(c)\simleq \no{\Box^ku}^2_{s+2l-\frac12}+\T{absorbable}+\T{neglectable},
\\
(a)+(d)\simleq \no{\Box^ku}^2_{s+2l-\frac12}+\T{absorbable}+\T{neglectable}.
\end{cases}
\end{equation*}
The first inequality is proved in the same way as the second part of 5.3 of \cite{BDKT06}. The second as in 5.4 of \cite{BDKT06} with the relevant change that
we do not have at our disposal their  estimate $|[\bar L,|z|^{2k}g_{1\bar 1}]|\simleq |z|^{2k-1-2(m-1)}$. Instead, we have to use, as a consequence of our key assumption \eqref{1.2}
\begin{equation*}
\begin{split}
[\bar L,|z|^{2k}g_{1\bar1}]&\simleq |z|^{2k-1}g_{1\bar1}+|z|^{2k}|g_{1\bar1\bar1}
\\
&\simleq |z|^{2k-1}g_{1\bar1}.
\end{split}
\end{equation*}
Thus, when we arrive at the two error terms in the second displayed formula of p. 692 (second terms in the third and fourth lines), we have the factor $z^{2k-1}g_{1\bar1}$. With the notations of our Lemma~\ref{l2.3}, we split this factor as $h=h_1h_2$ for $h_1=z^{2k-1}g_{1\bar1}^{\frac12}$ and $h_2=g_{1\bar1}^{\frac12}$ respectively and then control these error terms as sc (a) and lc (b).
The proof is complete.

\hskip15cm $\Box$

\end{document}